\documentclass[12pt,a4paper]{article}
\usepackage{amsfonts}
\usepackage{amssymb}
\usepackage{mathrsfs}
\usepackage{amsmath}
\usepackage{amsthm}
\usepackage{enumerate}
\usepackage{cite}
\usepackage[all]{xy}
\usepackage{extarrows}
\usepackage{indentfirst}
\allowdisplaybreaks
\newtheorem{defn}{Definition}[section]
\newtheorem{thm}[defn]{Theorem}
\newtheorem{lem}[defn]{Lemma}
\newtheorem{prop}[defn]{Proposition}
\newtheorem{cor}[defn]{Corollary}
\newtheorem{eg}[defn]{Example}
\newtheorem{re}[defn]{Remark}

\newcommand\relphantom[1]{\mathrel{\phantom{#1}}}
\newcommand{\bdefn}{\begin{defn}}
\newcommand{\edefn}{\end{defn}}
\newcommand{\bthm}{\begin{thm}}
\newcommand{\ethm}{\end{thm}}
\newcommand{\blem}{\begin{lem}}
\newcommand{\elem}{\end{lem}}
\newcommand{\bprop}{\begin{prop}}
\newcommand{\eprop}{\end{prop}}
\newcommand{\bcor}{\begin{cor}}
\newcommand{\ecor}{\end{cor}}
\newcommand{\beg}{\begin{eg}}
\newcommand{\eeg}{\end{eg}}
\newcommand{\bre}{\begin{re}}
\newcommand{\ere}{\end{re}}
\newcommand{\bpf}{\begin{proof}}
\newcommand{\epf}{\end{proof}}
\newcommand{\benu}{\begin{enumerate}}
\newcommand{\eenu}{\end{enumerate}}
\newcommand{\bc}{\begin{center}}
\newcommand{\ec}{\end{center}}
\newcommand{\bea}{\begin{eqnarray}}
\newcommand{\eea}{\end{eqnarray}}
\newcommand{\Bea}{\begin{eqnarray*}}
\newcommand{\Eea}{\end{eqnarray*}}
\newcommand{\beq}{\begin{equation}}
\newcommand{\eeq}{\end{equation}}
\newcommand{\Beq}{\begin{equation*}}
\newcommand{\Eeq}{\end{equation*}}
\newcommand{\bspl}{\begin{split}}
\newcommand{\espl}{\end{split}}

\begin{document}

\title{\textbf{On the structures of split Leibniz triple systems}
\author{ Yan Cao$^{1,2},$  Liangyun Chen$^{1}$
 \date{{\small {$^1$ School of Mathematics and Statistics, Northeast Normal
 University,\\
Changchun 130024, China}\\{\small {$^2$  Department of Basic
 Education,
 Harbin University of
Science and Technology,\\ Rongcheng Campus,  Rongcheng 264300,
China}}}}}} \maketitle
\date{}

\begin{abstract}

We study the structures of arbitrary split Leibniz triple systems.
By developing techniques of connections of roots for this kind of
triple systems,  under certain conditions, in the case of $T$ being
of maximal length, the simplicity
of the Leibniz triple systems is  characterized. \\

\end{abstract}
\renewcommand{\thefootnote}{\fnsymbol{footnote}}
\footnote[0]{ Corresponding author(L. Chen): chenly640@nenu.edu.cn.}
\footnote[0]{Supported by  NNSF of China (Nos. 11171055 and
11471090),  NSF of Jilin province (No. 201115006), Scientific
Research Fund of Heilongjiang Provincial Education Department
 (No. 12541184). }

\section{Introduction}

  Leibniz triple systems were introduced by Bremner and S\'{a}nchez-Ortega \cite{BS}. Leibniz triple systems were defined in a functorial manner using
 the Kolesnikov-Pozhidaev algorithm, which took the defining identities for a variety of algebras and produced the defining identities for the corresponding variety
  of dialgebras \cite{K}.  In \cite{BS}, Leibniz triple systems were obtained by applying the Kolesnikov-Pozhidaev algorithm to Lie triple systems. In \cite{BL5238999},
 Levi's theorem for Leibniz triple systems is determined.  Furthermore, Leibniz triple
  systems are related to Leibniz algebras in the same way that Lie triple systems related to Lie algebras. So it is natural to prove analogs of results from the theory
   of Lie triple systems to Leibniz triple systems.

In the present paper, we are interested in studying the structures of  arbitrary Leibniz triple systems  by focussing on the split ones. The class of the split ones is specially related to addition quantum
numbers, graded contractions, and deformations. For instance, for a physical system which displays
a symmetry of $T$, it is interesting to know in detail the structure of the split decomposition
because its roots can be seen as certain eigenvalues which are the additive quantum numbers
characterizing the state of such system. Recently, in \cite{BL52, BL52567, BL5234, BL528 }, the structures of arbitrary split Lie algebras, arbitrary split
  Leibniz algebras and arbitrary split  Lie triple systems have been determined by the techniques of connections of roots.
   Our work is essentially motivated by the work on  split Leibniz algebras and  split  Lie triple systems\cite{BL52567, BL5234}.

Throughout this paper, Leibniz triple systems $T$ are considered of
arbitrary dimension and over an arbitrary field $\mathbb{K}$. It is
worth to mention that, unless otherwise stated, there is not any
restriction on dim$T_{\alpha}$ or $\{k \in \mathbb{K}$:  $k \alpha
\in \Lambda^{1},$ for a fixed $\alpha \in \Lambda^{1}\}$, where
$T_{\alpha}$ denotes the root space associated to the root $\alpha$,
and $\Lambda^{1}$ the set of nonzero roots of $T$. This paper
proceeds as follows. In section 2, we establish the preliminaries on
split Leibniz triple systems theory. In section 3, we show that under certain conditions, in
the case of $T$ being of maximal length,  the simplicity
of the Leibniz triple systems is characterized.

\section{Preliminaries}

\bdefn{\rm\cite{BL5234}}  A \textbf{right Leibniz algebra} $L$ is a vector space over a field $\mathbb{K}$ endowed with a bilinear product
$[\cdot,\cdot]$ satisfying the Leibniz identity
$$[[y, z], x] = [[y, x], z] + [y, [z, x]],$$
for all $x, y, z \in L$.
\edefn

\bdefn{{\rm\cite{BS}}} A \textbf{Leibniz triple system} is a vector space $T$ endowed with a trilinear
operation $\{\cdot,\cdot,\cdot\}: T\times T\times T\rightarrow T$ satisfying
\begin{gather}
\{a,\{b,c,d\},e\}\!=\!\{\{a,b,c\},d,e\}
\!-\!\{\{a,c,b\},d,e\}\!-\!\{\{a,d,b\},c,e\}\!+\!\{\{a,d,c\},b,e\},\label{VIP1}\\
\{a,b,\{c,d,e\}\}\!=\!\{\{a,b,c\},d,e\}\!-\!\{\{a,b,d\},c,e\}\!-\!\{\{a,b,e\},c,d\}\!+\!\{\{a,b,e\},d,c\},\label{VIP2}
\end{gather}
for all $a, b, c, d, e \in T$.
\edefn

\beg\label{eg}
A Lie triple system gives a Leibniz triple system with the same ternary product.
If $L$ is a Leibniz algebra with product $[\cdot,\cdot]$, then $L$ becomes a Leibniz triple system by putting $\{x,y,z\}=[[x,y],z]$. More examples refer to \cite{BS}.
\eeg

\bdefn{{\rm\cite{BS}}} Let $I$ be a subspace of a Leibniz triple system  $T$. Then $I$ is called a \textbf{subsystem} of $T$, if $\{I,I,I\}\subseteq I;$ $I$ is
called an \textbf{ideal} of $T$,
if $\{I,T,T\}+\{T,I,T\}+\{T,T,I\}\subseteq I$.
\edefn

\bdefn
The \textbf{annihilator} of a Leibniz triple system $T$ is the set $\mathrm{Ann}(T)=\{x \in T: \{x, T, T\}+ \{T, x,  T\}+\{T,T,x\}=0\}$.
\edefn

\bprop{{\rm\cite{BL5238999}}}\label{38888}
Let $T$ be a  Leibniz triple system. Then the following assertions hold.

$\rm(1)$  $J$ is generated by $\{\{a,b,c\}-\{a,c,b\}+\{b,c,a\}: a,b,c \in T\}$, then $J$ is an ideal of $T$ satisfying $\{T,T,J\}=\{T,J,T\}=0$.

$\rm(2)$  $J$ is generated by $\{\{a,b,c\}-\{a,c,b\}+\{b,c,a\}: a,b,c \in T\}$, then $T$ is a Lie triple system if and only if $J=0$.

$\rm(3)$ $\{\{c,d,e\},b,a\}-\{\{c,d,e\},a,b\}-\{\{c,b,a\},d,e\}+\{\{c,a,b\},d,e\}-\{c,\{a,b,d\},e\}-\{c,d,\{a,b,e\}\}=0$, for all $a, b, c, d, e$ $\in$ $T$.
\eprop

\bdefn\label{ 3.16789} A Leibniz triple system $T$ is said to be
\textbf{simple} if its product is nonzreo and its only ideals are
$\{0\}$, $J$ and $T$, where $J$ is generated by
$\{\{a,b,c\}-\{a,c,b\}+\{b,c,a\}: a,b,c \in T\}$. \edefn

It should be noted that the above definition agrees with the
definition of a simple Lie triple system, since $J=\{0\}$ in this
case.

\bdefn{{\rm\cite{BS}}}\label{uni Leib envelop} The \textbf{standard embedding} of a Leibniz triple system $T$ is the two-graded right Leibniz algebra $L =
L^{0}\oplus L^{1}$, $L^{0}$ being the $\mathbb{K}$-$\rm span$ of $\{x \otimes y,  \ x, y \in T \}$, $L^{1}: =T$ and where the product is
given by
$$[(x \otimes y, z), (u \otimes v,w)]:= (\{x, y, u\} \otimes  v - \{x, y, v\} \otimes u
+ z \otimes w, \{x, y, w\} +\{z,u, v\}-\{z,v, u \}).$$
\edefn

Let us observe that $L^{0}$  with the product induced by the one in $L =
L^{0}\oplus L^{1}$  becomes a
 right Leibniz algebra.

\bdefn Let $T$ be a Leibniz triple system, $L = L^{0}\oplus L^{1}$ be its standard embedding, and $H^{0}$ be a  maximal abelian
subalgebra $($\rm MASA$)$
of $L^{0}$. For a linear function  $\alpha \in (H^{0})^{\ast},$
 we define the root space of $T$ $($with respect to
$H^{0}$$)$ associated to $\alpha$ as the subspace $T_{\alpha}:= \{t_{\alpha} \in T: [t_{\alpha}, h] = \alpha(h)t_{\alpha}$ for any h $\in H^{0}\}$.
The elements $\alpha \in (H^{0})^{\ast}$
 satisfying $T_{\alpha} \neq 0$ are called roots of $T$ with respect to $H^{0}$ and
we denote $\Lambda^{1}:= \{\alpha \in (H^{0})^{\ast}\setminus \{0\}: T_{\alpha}\neq 0 \}.$
\edefn

Let us observe that $T_{0} = \{t_{0} \in T: [t_{0}, h] = 0$ for any h $\in H^{0}\}$. In the following, we shall
denote by $\Lambda^{0}$ the set of all nonzero $\alpha \in (H^{0})^{\ast}$
such that $L_{\alpha}^{0}
:= \{ v_{\alpha}^{0}\in
 L^{0}: [v_{\alpha}^{0}, h]
 = \alpha(h)v_{\alpha}^{0}$
for any h $\in H^{0}\} \neq 0$.

 \blem \label{355} Let $T$ be a Leibniz triple system,  $L =L^{0}\oplus L^{1}$ be its standard embedding,
and  $H^{0}$ be a $\rm MASA$ of $L^{0}$. For $\alpha,\beta,\gamma
\in \Lambda^{1}\cup \{0\}$ and $\delta \in \Lambda^{0} \cup \{0\}$,
the following assertions hold.

$\rm(1)$ If $[T_{\alpha}, T_{\beta}]\neq 0$ then $\alpha+\beta \in  \Lambda^{0}\cup \{0\}$  and  $[T_{\alpha}, T_{\beta}]\subseteq L_{\alpha+\beta}^{0}$.

$\rm(2)$  If $[L_{\delta}^{0}, T_{\alpha}]\neq 0$ then $\delta+\alpha \in  \Lambda^{1}\cup \{0\}$  and  $[L_{\delta}^{0}, T_{\alpha}]\subseteq T_{\delta+\alpha}$.

$\rm(3)$  If $[T_{\alpha}, L_{\delta}^{0} ]\neq 0$ then $\alpha+\delta \in  \Lambda^{1}\cup \{0\}$  and  $[T_{\alpha}, L_{\delta}^{0}]\subseteq T_{\alpha+\delta}$.

$\rm(4)$ If $[L_{\delta}^{0}, L_{\gamma}^{0}]\neq 0$ then $\delta+\gamma \in  \Lambda^{0}\cup \{0\}$ and $[L_{\delta}^{0}, L_{\gamma}^{0} ]\subseteq L_{\delta+\gamma}^{0}$.

$\rm(5)$  If $\{T_{\alpha}, T_{\beta}, T_{\gamma}\}\neq 0$ then $\alpha+\beta+\gamma \in  \Lambda^{1}\cup \{0\}$  and  $\{T_{\alpha}, T_{\beta},T_{\gamma}\}\subseteq T_{\alpha+\beta+\gamma}$.
\elem

\bpf (1) For any $x \in T_{\alpha}$, $y \in T_{\beta}$ and $h \in H^{0}$, by Leibniz identity, one has $[[x,y],h]=[x,[y,h]]+[[x,h],y]=[x, \beta(h)y]+[\alpha(h)x,y]=(\alpha+\beta)(h)[x,y].$

(2)  For any $x \in L_{\delta}^{0}$, $y \in T_{\alpha}$ and $h \in H^{0}$, by Leibniz identity, one has $[[x,y],h]=[x,[y,h]]+[[x,h],y]=[x, \alpha(h)y]+[\delta(h)x,y]=(\delta+\alpha)(h)[x,y].$

(3)  For any $x \in T_{\alpha}$, $y \in L_{\delta}^{0}$,  and $h \in H^{0}$, by Leibniz identity, one has $[[x,y],h]=[x,[y,h]]+[[x,h],y]=[x, \delta(h)y]+[\alpha(h)x,y]=(\alpha+\delta)(h)[x,y].$

(4) For any $x \in L_{\delta}^{0}$, $y \in L_{\gamma}^{0}$ and $h \in H^{0}$, by Leibniz identity, one has $[[x,y],h]=[x,[y,h]]+[[x,h],y]=[x, \gamma(h)y]+[\delta(h)x,y]=(\delta+\gamma)(h)[x,y].$

(5) It is a consequence of Lemma $\ref{355}$ (1) and (2).
\epf

\bdefn \label{355468}Let $T$ be a Leibniz triple system, $L = L^{0}\oplus L^{1}$ be its standard embedding, and  $H^{0}$ be a $\rm MASA$
of $L^{0}$. We shall call that $T$ is a \textbf{split Leibniz triple system} $($with respect to $H^{0}$$)$ if $:$

$\rm(1)$ $T = T_{0}\oplus(\oplus_{\alpha \in \Lambda^{1}} T_{\alpha})$,

$\rm(2)$  $\{T_{0}, T_{0}, T_{0}\}=0$,

$\rm(3)$ $\{T_{\alpha}, T_{-\alpha}, T_{0}\}=0$, for $\alpha \in \Lambda^{1}$.

\noindent We say that $\Lambda^{1}$ is the root system of $T$.
\edefn

We also note that the facts $H^{0}\subset L^{0}=[T,T]$ and $T=T_{0}\oplus(\oplus_{\alpha \in \Lambda^{1}} T_{\alpha})$ imply
\beq\label{111}
H^{0}=[T_{0},T_{0}]+\sum_{\alpha \in \Lambda^{1}}[T_{\alpha},T_{-\alpha}].
\eeq

Finally, as $[T_{0}, T_{0}]\subset L^{0}_{0}=H^{0},$ we have
\beq\label{222} [T_{0},[T_{0},T_{0}]]=0. \eeq

We finally note that $\alpha \in \Lambda^{1}$ does not imply $\alpha \in \Lambda^{0}$.

\bdefn
A root system $\Lambda^{1}$ of a split Leibniz triple system $T$ is called \textbf{symmetric} if it satisfies that $\alpha \in  \Lambda^{1}$
implies $-\alpha \in \Lambda^{1}$.

A similar concept applies to the set $\Lambda^{0}$ of nonzero roots of $L^{0}$.
\edefn

In the following, $T$ denotes a split Leibniz triple system with a
symmetric root system $\Lambda^{1}$, and $T =
T_{0}\oplus(\oplus_{\alpha \in \Lambda^{1}} T_{\alpha})$ the
corresponding root decomposition.  We begin the study of split
Leibniz triple systems by developing the concept of connections of
roots.

\bdefn\label{333555}
Let $\alpha$ and $\beta$ be two nonzero roots, we shall say that $\alpha$ and $\beta$ are \textbf{connected} if there exists a family  $\{\alpha_{1},\alpha_{2},\cdots,\alpha_{2n},\alpha_{2n+1}\}\subset \Lambda^{1}\cup \{0\}$ of roots of $T$ such that

\noindent $\rm(1)$ $\{\alpha_{1},\alpha_{1}+\alpha_{2}+\alpha_{3},\alpha_{1}+\alpha_{2}+\alpha_{3}+\alpha_{4}+\alpha_{5},\cdots,\alpha_{1}+\cdots+\alpha_{2n}+\alpha_{2n+1}\} \subset \Lambda^{1},$

\noindent $\rm(2)$ $\{\alpha_{1}+\alpha_{2}, \alpha_{1}+\alpha_{2}+\alpha_{3}+\alpha_{4},\cdots, \alpha_{1}+\cdots+\alpha_{2n}\}\subset \Lambda^{0},$

\noindent $\rm(3)$  $\alpha_{1}=\alpha$ and $\alpha_{1}+\cdots+\alpha_{2n}+\alpha_{2n+1}= \pm\beta.$

\noindent We shall also say
that $\{\alpha_{1},\alpha_{2},\cdots,\alpha_{2n},\alpha_{2n+1}\}$ is a connection from $\alpha$ to $\beta$.
\edefn

We denote by
$$\Lambda_{\alpha}^{1}:= \{\beta \in \Lambda^{1}:  \alpha \ and \ \beta \ are \ connected\},$$
we can easily get that $\{\alpha\}$ is a connection from $\alpha$ to itself and to $-\alpha$. Therefore $\pm\alpha \in \Lambda_{\alpha}^{1}$.

\bdefn
A subset $\Omega^{1}$ of a root system $\Lambda^{1}$, associated to a split Leibniz triple system $T$, is called a\textbf{ root
subsystem} if it is symmetric, and for $\alpha, \beta, \gamma \in \Omega^{1}\cup \{0\} $ such that $\alpha + \beta \in \Lambda^{0}$ and
$\alpha +  \beta+ \gamma  \in  \Lambda^{1}$ then $\alpha +  \beta+ \gamma \in \Omega^{1}$.
\edefn

Let $\Omega^{1}$ be a root subsystem of $\Lambda^{1}$. We define
$$T_{0, \Omega^{1}} := \mathrm{span}_{\mathbb{K}}\{\{T_{\alpha}, T_{\beta}, T_{\gamma} \}:  \alpha +  \beta+ \gamma = 0;\  \alpha, \beta, \gamma \in  \Omega^{1}\cup\{0\}\} \subset T_{0}$$
and $V_{\Omega^{1}}:=\oplus_{\alpha \in \Omega^{1}}T_{\alpha}$. Taking into account the fact that $\{T_{0}, T_{0}, T_{0}\}=0$, it is straightforward to verify that
$T_{\Omega^{1}}:=T_{0,\Omega^{1}}\oplus V_{\Omega^{1}}$ is a  subsystem of $T$.  We will say that $ T_{\Omega^{1}}$ is a  subsystem associated to the
root subsystem $\Omega^{1}$.

\bprop \label{678}
If $\Lambda^{0}$ is symmetric, then the relation $\sim$ in $\Lambda^{1}$, defined by $\alpha \sim \beta$ if and only if $\beta \in \Lambda_{\alpha}^{1},$
is of equivalence.
\eprop

\bpf
This can be proved completely analogously to \cite[Proposition 3.1]{BL523}.
\epf

\bprop \label{67890}
 Let $\alpha$ be a nonzero root and suppose $\Lambda^{0}$ is symmetric. Then  $\Lambda_{\alpha}^{1}$
 is a root subsystem.
\eprop

\bpf This can be proved completely analogously to \cite[Lemma
3.1]{BL523}. \epf

\section{Split Leibniz triple system of maximal length. The simple case.}

In this section we focus on the simplicity of split Leibniz triple systems by centering our attention in those of maximal length. From now on char($\mathbb{K}$)=0.
\bdefn
We say that a split Leibniz triple systems $T$ is of \textbf{maximal length} if $\mathrm{dim}T_{\alpha}$=1 for any $\alpha \in \Lambda^{1}$.

\edefn

\blem \label{lemma 4.1}
Let $T$ be a split Leibniz triple systems. For any $\alpha, \beta$ $\in \Lambda^{1}$  with $\alpha \neq k\beta$, $k \in \mathbb{K}$, there
exists $h_{\alpha,\beta} \in H^{0}$ such that $\alpha(h_{\alpha,\beta}) \neq 0$ and $\beta( h_{\alpha,\beta}) = 0$.
\elem

\bpf As $\alpha \neq 0$, there exists $h \in H^{0} -\{0\}$ such that
$\alpha(h) \neq 0$. If $\beta(h) = 0$ we take $h_{\alpha,\beta}:=
h$. Suppose therefore $\beta(h)\neq 0$, let us write $k =
\frac{\alpha(h)}{\beta(h)}$. As $\alpha \neq k\beta$, there exists
$h^{'} \in H^{0}$  such that $\alpha(h^{'})\neq k\beta(h^{'}).$ we
can take $h_{\alpha,\beta}:=\beta(h^{'})h-\beta(h)h^{'}$. \epf

\blem \label{lemma 4.2}
Let $T = T_{0}\oplus(\oplus_{\alpha \in \Lambda^{1}} T_{\alpha})$ be a split Leibniz triple systems. If $I$ is an ideal of $T$ then $I=(I\cap T_{0})\oplus(\oplus_{\alpha \in \Lambda^{1}}(I\cap T_{\alpha})).$
\elem

\bpf
Let $x \in I$. We can write $x = t_{0}+\sum_{j=1}^{m}e_{\beta_{j}}$, for
$t_{0}\in T_{0}$, $e_{\beta_{j}}\in T_{\beta_{j}}$ and $\beta_{j} \neq \beta_{k}$ if $j\neq k$. Let us show that any $e_{\beta_{j}} \in I.$
 If $e_{\beta_{1}}=0$ then $e_{\beta_{1}} \in I$. Suppose  $e_{\beta_{1}}\neq 0$.
For any $\beta_{k_{r}}\neq p\beta_{1}, p\in \mathbb{K}$ and $k_{r}\in \{2,\cdots,m\}$, Lemma \ref{lemma 4.1} gives us $h_{\beta_{1},\beta_{k_{r}}} \in H^{0}$
satisfying  $\beta_{1}(h_{\beta_{1},\beta_{k_{r}}})\neq 0$ and $\beta_{k_{r}}(h_{\beta_{1},\beta_{k_{r}}})=0$.  From  here,

\begin{equation}\label{888}
[[\cdots[[x, h_{\beta_{1},\beta_{k_{2}}}],h_{\beta_{1},\beta_{k_{3}}}],\cdots],h_{\beta_{1},\beta_{k_{s}}}]=p_{1}e_{\beta_{1}}+\sum_{t=1}^{u}p_{k_{t}}e_{k_{t}\beta_{1}}\in I,
\end{equation}

\noindent where $p_{1}$, $k_{t}\in \mathbb{K}-\{0\}$, $k_{t}\neq 1$ and $p_{k_{t}} \in \mathbb{K}$.

If any $p_{k_{t}}=0,t=1,\cdots,u,$ then $p_{1}e_{\beta_{1}}\in I$ and so $e_{\beta_{1}}\in I$. Let us suppose some $p_{k_{t}}\neq 0$ and write (\ref{888}) as
\begin{equation}\label{8889}
p_{1}e_{\beta_{1}}+\sum_{t=1}^{v}p_{k_{t}}e_{k_{t}\beta_{1}}\in I,
\end{equation}
where $p_{1}, k_{t},p_{k_{t}}\in \mathbb{K}-\{0\}$, $k_{t} \neq 1, v\leq u$.

Let $h \in H^{0}$  such that $\beta_{1}(h)\neq 0$. Then

$$[  p_{1}e_{\beta_{1}}+\sum_{t=1}^{v}p_{k_{t}}e_{k_{t}\beta_{1}}, h]=p_{1}\beta_{1}(h)e_{\beta_{1}}+\sum_{t=1}^{v}p_{k_{t}}k_{t}\beta_{1}(h)e_{k_{t}\beta_{1}}\in I,$$

\noindent and so

\begin{equation}\label{888976}
p_{1}e_{\beta_{1}}+\sum_{t=1}^{v}p_{k_{t}}k_{t}e_{k_{t}\beta_{1}}\in I,\quad k_{t}\neq 1.
\end{equation}

 From (\ref{8889}) and (\ref{888976}), it follows easily that

\begin{equation}\label{88897655}
q_{1}e_{\beta_{1}}+\sum_{t=1}^{w}q_{k_{t}}e_{q_{t}\beta_{1}}\in I,
\end{equation}

 where $q_{1},q_{k_{t}} \in \mathbb{K}-\{0\}$, $q_{t} \in \{k_{t}: t=1,\cdots,v\}$ and $w<v$.

Following this process $($multiply (\ref{88897655})  with $h$ and compare the result with (\ref{88897655})  taking into account $q_{t}\neq 1$, and so
on$)$, we obtain $e_{\beta_{1}}\in I$.
 The same argument holds for any $\beta_{j},j\neq 1.$ From here, we deduce  $I=(I\cap T_{0})\oplus(\oplus_{\alpha \in \Lambda^{1}}(I\cap T_{\alpha})).$
\epf

Let us return to a split Leibniz triple system of maximal length $T$. From now on $T = T_{0}\oplus(\oplus_{\alpha \in \Lambda^{1}} T_{\alpha})$ denotes
 a split Leibniz triple system of maximal length. Using the previous Lemma, we assert that given any nonzero ideal $I$ of $T$ then

\begin{equation}\label{8889765566}
I=(I\cap T_{0})\oplus(\oplus_{\alpha \in \Lambda^{I}} T_{\alpha}),
\end{equation}
where $\Lambda^{I}:=\{\alpha \in  \Lambda^{1}: I\cap  T_{\alpha}\neq 0\}$.

In particular, case $I$=$J$, we get
\begin{equation}\label{888976556689}
J=(J\cap T_{0})\oplus(\oplus_{\alpha \in \Lambda^{J}} T_{\alpha}).
\end{equation}

\noindent From here, we can write
\begin{equation}\label{88897655668977}
\Lambda^{1}=\Lambda^{J}\cup \Lambda^{\neg J},
\end{equation}

\noindent where
$$\Lambda^{J}:=\{\alpha \in \Lambda^{1}: T_{\alpha}\subset J\}$$
\noindent and
$$\Lambda^{\neg J}:=\{\alpha \in \Lambda^{1}: T_{\alpha}\cap J=0\}.$$
 As consequence
\begin{equation}\label{8889765566897756}
T = T_{0}\oplus(\oplus_{\alpha \in \Lambda^{\neg J}} T_{\alpha})\oplus(\oplus_{\beta \in \Lambda^{J}} T_{\beta}).
\end{equation}

Next, we will consider $T$ satisfying  $\{T_{\alpha},T_{0},T_{\beta}\}\neq 0$, for $\alpha \in \Lambda^{J}$ and $\beta \in \Lambda^{\neg J}$. Under this assumption,
the fact that $T=\{T,T,T\}$, the split decomposition given by (\ref{8889765566897756}) and the Definition of a split Leibniz triple system  show

\begin{equation}\label{88897655668977565678}
T_{0} = \sum_{\alpha,\beta,\gamma \in \Lambda^{\neg J}\cup \{0\}\atop \alpha+\beta+\gamma=0}\{T_{\alpha},T_{\beta},T_{\gamma}\}.
\end{equation}

Now, observe that the concept of connectivity of nonzero roots given in Definition \ref{333555} is not strong
enough to detect if a given $\alpha \in \Lambda^{1}$ belongs to $\Lambda^{J}$ or to $\Lambda^{\neg J}$. Consequently we lose the information
respect to whether a given root space $T_{\alpha}$ is contained in $J$ or not, which is fundamental to study the
simplicity of $T$. So, we are going to refine the concept of connections of nonzero roots in the setup of
maximal length split Leibniz triple systems. The symmetry of $\Lambda^{J}$  and  $\Lambda^{\neg J}$ will be understood as usual. That
is, $\Lambda^{\gamma}$,$\gamma \in \{J, \neg J\}$, is called symmetric if $\alpha \in \Lambda^{\gamma}$ implies $-\alpha \in \Lambda^{\gamma}$.

In the following, T denotes a split Leibniz triple system whose root space satisfies $\{T_{\alpha},T_{0},T_{\beta}\}\neq 0$, where $\alpha \in \Lambda^{J}$ and $\beta \in \Lambda^{\neg J}$.

\bdefn
Let $\alpha, \beta \in \Lambda^{\gamma}$  with $\gamma  \in \{J, \neg J\}$. We say that $\alpha$ is $\neg J$-connected to $\beta$, denoted by
$\alpha\sim _{\neg J} \beta$, if there exist
$$\alpha_{2}, \cdots, \alpha_{2n+1} \in \Lambda^{\neg J}\cup \{0\}$$
\noindent such that

\noindent $(\rm 1)$ $\{\alpha_{1},\alpha_{1}+\alpha_{2}+\alpha_{3},\alpha_{1}+\alpha_{2}+\alpha_{3}+\alpha_{4}+\alpha_{5},\cdots,\alpha_{1}+\cdots\alpha_{2n}+\alpha_{2n+1}\} \subset \Lambda^{\gamma}, $

\noindent $(\rm 2)$ $\{\alpha_{1}+\alpha_{2}, \alpha_{1}+\alpha_{2}+\alpha_{3}+\alpha_{4},\cdots, \alpha_{1}+\cdots\alpha_{2n}\}\subset \Lambda^{0},$

\noindent $(\rm 3)$ $\alpha_{1}=\alpha$ and $\alpha_{1}+\cdots+\alpha_{2n}+\alpha_{2n+1}= \pm\beta.$

 We shall also say
that $\{\alpha_{1},\alpha_{2},\cdots,\alpha_{2n},\alpha_{2n+1}\}$ is a $\neg J$-connection from $\alpha$ to $\beta$.

\edefn

\bprop \label{67890345777777}
 The following assertions  hold.

$(\rm 1)$ If  $\Lambda^{\neg J}$  is symmetric, then the relation $\sim_{\neg J}$ is an equivalence relation in $\Lambda^{\neg J}$.

$(\rm 2)$ If $ T= \{T,T,T\}$ and $\Lambda^{\neg J}$, $\Lambda^{J}$ are symmetric, then the relation $\sim_{\neg J}$ is an equivalence relation in $\Lambda^{J}$.
\eprop

\bpf
(1) Can be proved in a similar way to Proposition \ref{678}.

(2)   Note that $T_{0} = \sum_{\alpha,\beta,\gamma \in \Lambda^{\neg J}\cup \{0\}\atop \alpha+\beta+\gamma=0}\{T_{\alpha},T_{\beta},T_{\gamma}\}$ and $H^{0}=[T_{0}, T_{0}]+\sum_{\alpha \in \Lambda^{\neg J}}[T_{\alpha}, T_{-\alpha}]+\sum_{\alpha \in \Lambda^{ J}}[T_{\alpha}, T_{-\alpha}]$. Let $\delta \in \Lambda^{J}$, one gets
$$[T_{\delta},\sum_{\alpha \in \Lambda^{ J}}[T_{\alpha},T_{-\alpha}]]\subset \sum_{\alpha \in \Lambda^{J}}[[T_{\delta},  T_{\alpha}],T_{-\alpha}]+ \sum_{\alpha \in \Lambda^{J}}[[T_{\delta},  T_{-\alpha}],T_{\alpha}] = 0.$$
Since $\delta \neq 0$,
one gets either $[T_{\delta},[T_{0},T_{0}]]\neq 0$ or $[T_{\delta},\sum_{\alpha \in \Lambda^{\neg J}}[T_{\alpha},T_{-\alpha}]]\neq 0$.

Suppose  $[T_{\delta},[T_{0},T_{0}]]\neq 0$. By Leibniz identity, one gets $[[T_{\delta},T_{0}],T_{0}]\neq 0.$ Then $\neg J$-connection $\{\delta,0,0\}$ gives us $\delta\sim_{\neg J} \delta$.

Suppose  $[T_{\delta},\sum_{\alpha \in \Lambda^{\neg J}}[T_{\alpha},T_{-\alpha}]]\neq 0$. By Leibniz identity, either $[[T_{\delta},T_{\alpha}],T_{-\alpha}]\neq 0$ or
$[[T_{\delta},T_{-\alpha}],T_{\alpha}]\neq 0.$ In the first case, the $\neg J$-connection $\{\delta, \alpha, -\alpha\}$ gives us $\delta\sim_{\neg J} \delta$ while in the second case the  $\neg J$-connection $\{\delta, -\alpha, \alpha\}$ gives us the same result.

Consequently,  $\sim_{\neg J}$ is reflexive in $\Lambda^{J}$.  The symmetric and transitive character
of $\sim_{\neg J}$  in $\Lambda^{J}$ follows as in Proposition \ref{678}.
\epf

Let us introduce the notion of root-multiplicativity in the framework of split Leibniz triple systems of
maximal length, in a similar way to the ones for split Leibniz algebras and split Lie triple systems $($see
\cite {BL5234,BL52567} for these notions and examples$)$.

\bdefn
We say that a split Leibniz  triple system of maximal length $T$ is root-multiplicative if the
below conditions hold.

$(\rm 1)$ Given $\alpha, \beta, \gamma \in \Lambda^{\neg J}\cup \{0\}$ such that $\alpha+\beta \in \Lambda^{0}$ and  $\alpha+\beta+\gamma \in \Lambda^{1},$ then $\{T_{\alpha}, T_{\beta}, T_{\gamma}\}\neq 0$.

$(\rm 2)$ Given  $\alpha, \beta \in \Lambda^{\neg J}\cup \{0\}$ and $\gamma \in \Lambda^{J}$ such that  $\alpha+\beta \in \Lambda^{0}$ and  $\alpha+\beta+\gamma \in \Lambda^{J}$, then $\{T_{\gamma}, T_{\beta}, T_{\alpha}\}\neq 0$.
\edefn

\blem \label{lemma 4.4}
Let $T$ be a root-multiplicative split Leibniz triple systems with $\mathrm{Ann} (T)=0$. If for any $\alpha\in \Lambda^{1}$, we have $ \mathrm{dim}L_{\alpha}^{0}=1 $. Then there is not any nonzero ideal of $T$ contained in $T_{0}$.
\elem
\bpf
Suppose there exists a nonzero ideal of $T$ such that $I\subset T_{0}$. The facts that $\{T_{0}, T_{0}, T_{0}\}=0$ gives $\{I,  T_{0}, T_{0}\}=0$, $\{ T_{0}, I, T_{0}\}=0$ and $\{ T_{0}, T_{0},I \}=0$. Given $\alpha \in \Lambda^{1}$, since $$\{I, T_{0},  T_{\alpha}\}+\{I,  T_{\alpha},  T_{0}\}\subset  T_{\alpha}\cap T_{0}=0,$$
$$\{ T_{0}, I,  T_{\alpha}\}+\{  T_{\alpha},  I, T_{0}\}\subset  T_{\alpha}\cap T_{0}=0,$$
 and
 $$\{ T_{0},   T_{\alpha}, I\}+\{  T_{\alpha},  T_{0}, I\}\subset  T_{\alpha}\cap T_{0}=0,$$
one gets
$$\{I, T_{0},  T_{\alpha}\}=\{I,  T_{\alpha},  T_{0}\}=0,$$
$$\{ T_{0}, I,  T_{\alpha}\}=\{  T_{\alpha}, I, T_{0}\}=0,$$
and
$$\{ T_{0},  T_{\alpha} ,I\}=\{ T_{\alpha},  T_{0}, I \}=0.$$
Given also $\beta \in \Lambda^{1}$, with $\alpha+\beta\neq 0$, one has
$$\{I, T_{\alpha},  T_{\beta}\}\subset    T_{\alpha+\beta}\cap  T_{0}=0,$$
$$\{ T_{\alpha}, I, T_{\beta}\}\subset    T_{\alpha+\beta}\cap  T_{0}=0,$$
$$\{ T_{\alpha}, T_{\beta}, I\}\subset  T_{\alpha+\beta}\cap  T_{0}=0.$$
We also have $\{ T_{\alpha},  T_{-\alpha}, I\}\subset\{ T_{\alpha},  T_{-\alpha}, T_{0}\}=0$ $($ see the Definition of a split Leibniz triple system$)$.
As $\mathrm{Ann}(T)=0$, one gets
$$\{I, T_{\alpha},  T_{-\alpha}\}\neq 0$$
 or
$$\{ T_{\alpha}, I,  T_{-\alpha}\}\neq 0.$$
We treat separately two cases.

Case 1:  $\{I, T_{\alpha},  T_{-\alpha}\}\neq 0$. Thus, there exist $t_{\pm \alpha}\in T_{\pm \alpha}$ and $t_{0}\in I$ such that $\{t_{0}, t_{\alpha}, t_{-\alpha}\}\neq 0$. Hence, $0 \neq [{t_{0}, t_{\alpha}]\in L_{\alpha}^{0}}$.
  Using dim$L_{\alpha}^{0}=1$ and the root-multiplicativity  of $T$ $($consider the roots 0, $\alpha$, $0\in \Lambda^{1}\cup \{0\}$$)$, there exists $0 \neq  t_{0}^{'} \in T_{0}$ such that $0 \neq  \{t_{0}, t_{\alpha},  t_{0}^{'}\}\in T_{\alpha}$. As $t_{0}\in I$, we conclude $0 \neq  t_{\alpha}^{'}:=\{t_{0}, t_{\alpha},  t_{0}^{'}\}\in I\subset T_{0}$, a contradiction. Hence $I$ is not contained in $T_{0}$.

Case 2:  $\{ T_{\alpha}, I,  T_{-\alpha}\}\neq 0$. Thus, there exist $t_{\pm \alpha}\in T_{\pm \alpha}$ and $t_{0}\in I$ such that $\{ t_{\alpha}, t_{0}, t_{-\alpha}\}\neq 0$. Hence, $0 \neq [{t_{\alpha}, t_{0}]\in L_{\alpha}^{0}}$.
  Using dim$L_{\alpha}^{0}=1$ and the root-multiplicativity  of $T$ $($consider the roots  $\alpha$, 0, $0\in \Lambda^{1}\cup \{0\}$$)$, there exists $0 \neq  t_{0}^{'} \in T_{0}$ such that $0 \neq  \{t_{\alpha}, t_{0},   t_{0}^{'}\}\in T_{\alpha}$. As $t_{0}\in I$, we conclude $0 \neq  t_{\alpha}^{'}:=\{t_{\alpha}, t_{0},   t_{0}^{'}\}\in I\subset T_{0}$, a contradiction. Hence $I$ is not contained in $T_{0}$.
\epf

 Another interesting notion related to split Leibniz triple
systems of maximal length $T$ is those of Lie-annihilator. Write $T
= T_{0}\oplus(\oplus_{\alpha \in \Lambda^{\neg J}}
T_{\alpha})\oplus(\oplus_{\beta \in \Lambda^{J}} T_{\beta})$ (see
(\ref{8889765566897756})).

\bdefn
The Lie-annihilator of a split Leibniz triple system of maximal length $T$ is the set
\begin{align*}
\mathrm{Ann}_{\mathrm{Lie}}(T)=&\Big\{x \in T: \{x, T_{0}\oplus(\oplus_{\alpha \in \Lambda^{\neg J}} T_{\alpha}),T_{0}\oplus(\oplus_{\alpha \in \Lambda^{\neg J}} T_{\alpha})\}\\
+&\{T_{0}\oplus(\oplus_{\alpha \in \Lambda^{\neg J}} T_{\alpha}), x, T_{0}\oplus(\oplus_{\alpha \in \Lambda^{\neg J}} T_{\alpha})\}\\
+&\{T_{0}\oplus(\oplus_{\alpha \in \Lambda^{\neg J}} T_{\alpha}), T_{0}\oplus(\oplus_{\alpha \in \Lambda^{\neg J}} T_{\alpha}), x\}=0\Big\}.
\end{align*}
\edefn

Clearly the above Definition agrees with the Definition of annihilator of a Lie triple system, since in this
case $\Lambda^{J}=\emptyset$. We also have $\mathrm{Ann}(T)\subset \mathrm{Ann}_{\mathrm{Lie}}(T)$.

\bprop \label{67890345777777999}
Suppose $T= \{T, T,T\}$ and  $T$ is root-multiplicative. If $\Lambda^{\neg J}$
 has all of its roots $\neg J$-connected,
then any ideal $I$ of $T$ such that $I \not \subseteq T_{0}\oplus J$ satisfies $I = T$.
\eprop

\bpf
By (\ref{8889765566}) and (\ref{88897655668977}), we can write
$$I=(I\cap T_{0})\oplus(\oplus_{\alpha_{i} \in \Lambda^{\neg J,I}} T_{\alpha_{i}})\oplus(\oplus_{\beta_{j} \in \Lambda^{ J,I}} T_{\beta_{j}}),$$
\noindent where $\Lambda^{\neg J,I}:=\Lambda^{\neg J}\cap \Lambda^{I}$ and $\Lambda^{ J,I}:=\Lambda^{ J}\cap \Lambda^{I}.$  As $I \not \subseteq T_{0}\oplus J$, one gets $
\Lambda^{\neg J,I}\neq \emptyset$ and so we can fix some $\gamma_{0} \in \Lambda^{\neg J,I}$ such that
\begin{equation}\label{888976444}
T_{\gamma_{0}}\subset I.
\end{equation}

For any $\beta \in \Lambda^{\neg J}$, $\beta\neq \pm \gamma_{0}$. The fact that $\gamma_{0}$ and $\beta$ are $\neg J$-connected gives us a  $\neg J$-connection
$\{\gamma_{1},\cdots,\gamma_{2r+1}\}\subset \Lambda^{\neg J}\cup \{0\}$ from $\gamma_{0}$ to $\beta$ such that

$\gamma_{1}=\gamma_{0},$

$\gamma_{1}+\gamma_{2}, \gamma_{1}+\gamma_{2}+\gamma_{3}+\gamma_{4},\cdots, \gamma_{1}+\cdots+\gamma_{2r}\in \Lambda^{0},$

 $\gamma_{1},\gamma_{1}+\gamma_{2}+\gamma_{3},\cdots,\gamma_{1}+\cdots+\gamma_{2r}+\gamma_{2r+1} \in \Lambda^{\neg J}$,

\noindent and $\gamma_{1}+\cdots+\gamma_{2r}+\gamma_{2r+1} = \pm\beta$.

Consider $\gamma_{0}=\gamma_{1}, \gamma_{2},\gamma_{3}$ and $\gamma_{1}+\gamma_{2}+\gamma_{3}$. Since $\gamma_{1}, \gamma_{2},\gamma_{3} \in \Lambda^{\neg J}\cup \{0\},$ the root-multiplicativity and maximal length
of $T$ show $\{T_{\gamma_{0}}, T_{\gamma_{2}}, T_{\gamma_{3}}\}=T_{\gamma_{0}+\gamma_{2}+\gamma_{3}},$ and by  (\ref{888976444}),
$$T_{\gamma_{0}+\gamma_{2}+\gamma_{3}}\subset I.$$
\noindent Following this process with the $\neg J$-connection
$\{\gamma_{1},\cdots,\gamma_{2r+1}\}$,  we obtain that
$$T_{\gamma_{0}+\gamma_{2}+\gamma_{3}+\cdots+\gamma_{2r+1}}\subset I.$$
\noindent From here, we get that either
\begin{equation}\label{888976444469}
 T_{\beta}\subset I\ \mathrm{or} \ T_{-\beta}\subset I,
\end{equation}
for any $\beta \in \Lambda^{\neg J}$.

Observe that as a consequence of $T= \{T, T,T\}$ and  $\{T_{\alpha},T_{0},T_{\beta}\}\neq 0$  where $\alpha \in \Lambda^{J}$, $\beta \in \Lambda^{\neg J},$ one gets
\begin{equation}\label{888976444347}
T_{0} = \sum_{\alpha,\beta,\gamma \in \Lambda^{\neg J}\cup \{0\}\atop \alpha+\beta+\gamma=0}\{T_{\alpha},T_{\beta},T_{\gamma}\}.
\end{equation}

Let us study the products $\{T_{\alpha},T_{\beta},T_{\gamma}\}$ of (\ref{888976444347}) in order to show $T_{0}\subset I$. Taking into
account  Definition \ref{355468} (2), (3) and the fact that $\alpha+\beta+\gamma=0$  with $\alpha,\beta,\gamma \in \Lambda^{\neg J}\cup \{0\},$ we
can suppose $\gamma\neq 0 $ and either $\alpha \neq  0$ or $\beta \neq 0$. Suppose  $\alpha \neq  0$ and $\beta = 0$ (resp. $\alpha = 0$ and
 $\beta \neq 0$), one gets $\alpha = -\gamma$ (resp. $\beta =-\gamma$), and by  (\ref{888976444469}), $\{T_{\alpha},T_{\beta},T_{\gamma}\}$  =$\{T_{-\gamma},T_{0},T_{\gamma}\} \subset I$,
(resp. $\{T_{\alpha},T_{\beta},T_{\gamma}\}$  =$\{T_{0},T_{-\gamma},T_{\gamma}\} \subset I$). If the three elements in $\{\alpha, \beta, \gamma\}$ are nonzero, in
case some $T_{\epsilon}\subset I $,  $\epsilon \in  \{\alpha, \beta, \gamma\}$, then clearly  $\{T_{\alpha},T_{\beta},T_{\gamma}\}\subset I.$

 Finally, consider the case
in which any of the $T_{\epsilon}$  does not belong to $I$. If  $\{T_{\alpha},T_{\beta},T_{\gamma}\} = 0$ then $\{T_{\alpha},T_{\beta},T_{\gamma}\}\subset I $.
If $\{T_{\alpha},T_{\beta},T_{\gamma}\}\neq  0$, necessarily $\alpha + \beta \neq 0$ and so $\alpha + \beta \in \Lambda^{0}$. From here,  by
root-multiplicativity, one gets $0 \neq \{T_{\alpha},T_{\beta},T_{-\beta}\} =T_{\alpha}$. By  (\ref{888976444469}), we have $T_{\alpha} \subset I$
and so $\{T_{\alpha},T_{\beta},T_{\gamma}\} \subset I$. Therefore (\ref{888976444347}) implies
\begin{equation}\label{8889764443475778}
T_{0}\subset I.
\end{equation}

Now, given any $\delta \in \Lambda^{1}$, the facts $\delta \neq 0$, $T_{0}\subset I$, root-multiplicativity and $H^{0}=[T_{0}, T_{0}]+\sum_{\alpha \in \Lambda^{\neg J}}[T_{\alpha}, T_{-\alpha}]+\sum_{\alpha \in \Lambda^{ J}}[T_{\alpha}, T_{-\alpha}]$ show  either $[T_{\delta},[T_{0},T_{0}]]\neq 0$ or $[T_{\delta},[T_{\alpha},T_{-\alpha}]]\neq 0,$ for $\alpha \in \Lambda^{\neg J}$. We will distinguish respectively.

Suppose $[T_{\delta},[T_{0},T_{0}]]\neq 0.$ By Leibniz identity, we have $[T_{\delta},[T_{0},T_{0}]]\subseteq [[T_{\delta},T_{0}],T_{0}].$ Since $T_{0}\subset I$, one gets $[[T_{\delta},T_{0}],T_{0}]\subset I$, and so  $[T_{\delta},[T_{0},T_{0}]]\subset I.$ By  root-multiplicativity, $0 \neq [T_{\delta},[T_{0},T_{0}]]=T_{\delta}$, so
\begin{equation}\label{8889764443475778555}
T_{\delta}\subset I,
\end{equation}
for $\delta \in \Lambda^{1}$.

Suppose $[T_{\delta},[T_{\alpha},T_{-\alpha}]]\neq 0$, for $\alpha \in \Lambda^{\neg J}.$  By Leibniz identity, we have $[T_{\delta},[T_{\alpha},T_{-\alpha}]]\subseteq [[T_{\delta},T_{\alpha}],T_{-\alpha}]+[[T_{\delta},T_{-\alpha}],T_{\alpha}].$ By (\ref{888976444469}), one gets $[[T_{\delta},T_{\alpha}],T_{-\alpha}]\subset I$ and $[[T_{\delta},T_{-\alpha}],T_{\alpha}]\subset I$ for  $\alpha \in \Lambda^{\neg J}.$  By  root-multiplicativity, one gets $0 \neq [T_{\delta},[T_{\alpha},T_{-\alpha}]]=T_{\delta}$, and so
\begin{equation}\label{8885809764443475778555}
T_{\delta}\subset I,
\end{equation}
 for $\delta \in \Lambda^{1}$.

 From (\ref{888976444}),  (\ref{8889764443475778555}) and (\ref{8885809764443475778555}), we conclude $I=T$.
\epf

\bprop \label{67891110345777777999}
Suppose $T = \{T,T,T\}$, $\mathrm{Ann}(T) = 0$ and $T$ is root-multiplicative. If  $\Lambda^{\neg J}$, $\Lambda^{J}$ are symmetric,
 $\Lambda^{ J}$ has all of its roots $\neg J$-connected and $\{T_{0}, T_{\alpha}, T_{\beta}\}=0,$ for $\alpha, \beta \in  \Lambda^{\neg J}$, then any nonzero ideal $I$ of $T$ such that $I \subseteq J$ satisfies either
$I = J$ or $J = I \oplus K$ with $K$ an ideal of $T$.
\eprop

\bpf
By (\ref{8889765566}) and (\ref{88897655668977}), we can write
$$I=(I\cap T_{0})\oplus(\oplus_{\alpha_{i} \in \Lambda^{ J,I}} T_{\alpha_{i}}),$$
where $\Lambda^{ J,I}\subset \Lambda^{ J}$. Observe that the fact $\mathrm{Ann}(T) = 0$ implies
\begin{equation}\label{888778555}
J\cap T_{0}=0.
\end{equation}
Indeed, from the fact that $\{ T_{0}, T_{0}, T_{0}\}=0$, one gets
\begin{equation}\label{888778555777}
\{J\cap T_{0}, T_{0}, T_{0}\}+\{T_{0}, J\cap T_{0},  T_{0}\}+\{ T_{0}, T_{0},J\cap T_{0}\}=0.
\end{equation}
By Proposition \ref{38888} (1), it is easy to see
\begin{equation}\label{88877811777}
\{J\cap T_{0}, T_{0}, \oplus_{\alpha \in \Lambda^{J}}T_{\alpha}\}+\{T_{0}, J\cap T_{0},  \oplus_{\alpha \in \Lambda^{J}}T_{\alpha}\}+\{T_{0}, \oplus_{\alpha \in \Lambda^{J}} T_{\alpha}, J\cap T_{0}\}=0.
\end{equation}
By Proposition \ref{38888} (1), it is easy to see
 $$\{T_{0}, J\cap T_{0},  \oplus_{\beta \in \Lambda^{\neg J}}T_{\beta}\}=0$$
 and
$$\{T_{0},  \oplus_{\beta \in \Lambda^{\neg J}}T_{\beta}, J\cap T_{0} \}=0.$$

\noindent As $\{J\cap T_{0}, T_{0}, \oplus_{\beta \in \Lambda^{\neg J}}T_{\beta}\}\subset  \oplus_{\beta \in \Lambda^{\neg J}}T_{\beta}\cap J=0,$
  one gets
  $$\{J\cap T_{0}, T_{0}, \oplus_{\beta \in \Lambda^{\neg J}}T_{\beta}\}=0.$$
  
 Therefore, one gets
\begin{equation}\label{88875800078475611777}
\{J\cap T_{0}, T_{0}, \oplus_{\beta \in \Lambda^{\neg J}}T_{\beta}\}+\{T_{0}, J\cap T_{0},  \oplus_{\beta \in \Lambda^{\neg J}}T_{\beta}\}+\{T_{0},  \oplus_{\beta \in \Lambda^{\neg J}}T_{\beta}, J\cap T_{0} \}=0.
\end{equation}
By Proposition \ref{38888} (1), it is easy to see
\begin{equation}\label{8887781177778000}
\{J\cap T_{0},  \oplus_{\alpha \in \Lambda^{J}}T_{\alpha} ,T_{0}\}+\{\oplus_{\alpha \in \Lambda^{J}}T_{\alpha}, J\cap T_{0} , T_{0}\}+\{\oplus_{\alpha \in \Lambda^{J}} T_{\alpha}, T_{0},  J\cap T_{0}\}=0.
\end{equation}
Similarly, we also get
 \begin{equation}\label{8887781177777776668000}
\{J\cap T_{0},  \oplus_{\alpha \in \Lambda^{J}}T_{\alpha}, \oplus_{\beta\in \Lambda^{J}}T_{\beta}\}=0,
\end{equation}

\begin{equation}\label{88}
\{\oplus_{\alpha \in \Lambda^{J}}T_{\alpha}, J\cap T_{0} , \oplus_{\beta\in \Lambda^{J}}T_{\beta}\}=0,
\end{equation}

\begin{equation}\label{80}
\{\oplus_{\alpha \in \Lambda^{J}} T_{\alpha}, \oplus_{\beta\in \Lambda^{J}}T_{\beta},  J\cap T_{0}\}=0,
\end{equation}

\begin{equation}\label{888744478666481177776668000}
\{J\cap T_{0},  \oplus_{\alpha \in \Lambda^{J}}T_{\alpha}, \oplus_{\beta\in \Lambda^{\neg J}}T_{\beta}\}=0,
\end{equation}

\begin{equation}\label{8000}
\{\oplus_{\alpha \in \Lambda^{J}}T_{\alpha}, J\cap T_{0} , \oplus_{\beta\in \Lambda^{\neg J}}T_{\beta}\}=0,
\end{equation}

\begin{equation}\label{888800}
\{\oplus_{\alpha \in \Lambda^{J}} T_{\alpha}, \oplus_{\beta\in \Lambda^{\neg J}}T_{\beta},  J\cap T_{0}\}=0,
\end{equation}

\begin{equation}\label{33378000}
\{J\cap T_{0},  \oplus_{\beta \in \Lambda^{\neg J}}T_{\beta} ,T_{0}\}+\{\oplus_{\beta \in \Lambda^{\neg J}}T_{\beta}, J\cap T_{0}, T_{0}\}+\{\oplus_{\beta \in \Lambda^{\neg J}}T_{\beta}, T_{0},  J\cap T_{0}\}=0,
\end{equation}

\begin{equation}\label{33378000369}
\{J\cap T_{0},   \oplus_{\beta\in \Lambda^{\neg J}}T_{\beta}, \oplus_{\alpha \in \Lambda^{J}}T_{\alpha}\}=0,
\end{equation}

\begin{equation}\label{0369}
\{\oplus_{\beta\in \Lambda^{\neg J}}T_{\beta}, J\cap T_{0}, \oplus_{\alpha \in \Lambda^{J}}T_{\alpha}\}=0,
\end{equation}

\begin{equation}\label{339}
\{\oplus_{\beta\in \Lambda^{\neg J}T_{\beta}}, \oplus_{\alpha \in \Lambda^{J}}T_{\alpha},   J\cap T_{0}\}=0.
\end{equation}

\noindent By Proposition \ref{38888} (1), it is easy to see
\begin{equation}\label{5679}
 \{\oplus_{\beta\in \Lambda^{\neg J}}T_{\beta}, J\cap T_{0}, \oplus_{\alpha \in \Lambda^{\neg J}}T_{\alpha}\}=0,
 \end{equation}

\noindent and
\begin{equation}\label{1179}
 \{\oplus_{\beta\in \Lambda^{\neg J}}T_{\beta}, \oplus_{\alpha \in \Lambda^{\neg J}}T_{\alpha}, J\cap T_{0} \}=0.
 \end{equation}

\noindent By known condition, we have
\begin{equation}\label{118889}
 \{J\cap T_{0}, \oplus_{\beta\in \Lambda^{\neg J}}T_{\beta}, \oplus_{\alpha \in \Lambda^{\neg J}}T_{\alpha}\}=0.
  \end{equation}

\noindent From  (\ref{888778555777}), (\ref{88877811777}), (\ref{88875800078475611777}), (\ref{8887781177778000}), (\ref{8887781177777776668000}), (\ref{88}), (\ref{80}), (\ref{888744478666481177776668000}), (\ref{8000}), (\ref{888800}),  (\ref{33378000}), (\ref{33378000}),  (\ref{33378000369}), (\ref{0369}), (\ref{339}), (\ref{5679}), (\ref{1179}), (\ref{118889}) and (\ref{8889765566897756}), one gets
 $$\{J\cap T_{0},T,T\}+\{T,J\cap T_{0},T\}+\{T,T,J\cap T_{0}\}=0.$$
From here $J\cap T_{0}\subset$ $\mathrm{Ann}(T)=0$. Hence, we can write
 $$I=\oplus_{\alpha_{i} \in \Lambda^{ J,I}} T_{\alpha_{i}},$$
 with $\Lambda^{ J,I}\neq \emptyset$, and so we can take some $\alpha_{0}\in \Lambda^{ J,I}$ such that  $T_{\alpha_{0}}\subset I$. We can argue with the
  root-multiplicativity
and the maximal length of $T$ as in Proposition \ref{67890345777777999} to conclude that given any $\beta \in \Lambda^{ J}$,
there exists a $\neg J$-connection $\{\gamma_{1},\cdots, \gamma_{2r+1}\}$ from $\alpha_{0}$ to $\beta$ such that
$$\{\{\cdots,\{\{T_{\alpha_{0}}, T_{\gamma_{2}}, T_{\gamma_{3}}\},T_{\gamma_{4}},T_{\gamma_{5}}\},\cdots\},T_{\gamma_{2r}}, T_{\gamma_{2r+1}}\}= T_{\pm \beta}$$
and so $T_{\epsilon \beta}\subset I$ for some $\epsilon \in \pm 1$. That is
\begin{equation}\label{3347009}
\epsilon_{\beta}\beta \in \Lambda^{ J,I} \ \mathrm{for} \ \mathrm{any} \ \beta \in \Lambda^{J}  \ \mathrm{and} \ \mathrm{some} \ \epsilon_{\beta}\in \pm 1.
  \end{equation}

 Suppose  $-\alpha_{0}\in \Lambda^{ J,I}$. Then we also have that $\{-\gamma_{1},\cdots, -\gamma_{2r+1}\}$ from $\alpha_{0}$ to $\beta$  is a $\neg J$-connection from
 $-\alpha_{0}$ to $\beta$ satisfying
 $$\{\{\cdots,\{\{T_{-\alpha_{0}}, T_{-\gamma_{2}}, T_{-\gamma_{3}}\},T_{-\gamma_{4}},T_{-\gamma_{5}}\},\cdots\},T_{-\gamma_{2r}}, T_{-\gamma_{2r+1}}\}= T_{-\epsilon_{\beta}\beta}\subset I$$
 and so $T_{\beta}+T_{-\beta}\subset I$. Hence, (\ref{888976556689}) and  (\ref{888778555}) imply that $I=J$.

 Now,  suppose there is not any  $\alpha_{0}\in \Lambda^{ J,I}$ such that  $-\alpha_{0}\in \Lambda^{ J,I}$.  (\ref{3347009}) allows us to write $\Lambda^{ J}=
 \Lambda^{ J,I}\cup (-\Lambda^{ J,I})$ and (together with (\ref{888976556689}) and  (\ref{888778555})) asserts that by denoting $K = \oplus_{\alpha_{i} \in \Lambda^{ J,I}} T_{-\alpha_{i}}$, we have
$$J = I \oplus K.$$

 Let us finally show that $K$ is an ideal of $T$. We have $\{T,K,T\}+\{T,T,K\}=0$ and
 \begin{align*}
 &\{K,T,T\}\\
 =&\{K,  T_{0}\oplus(\oplus_{\alpha \in \Lambda^{\neg J}} T_{\alpha})\oplus(\oplus_{\beta \in \Lambda^{J}} T_{\beta}), T_{0}\oplus(\oplus_{\alpha \in \Lambda^{\neg J}} T_{\alpha})\oplus(\oplus_{\beta \in \Lambda^{J}} T_{\beta})\}\\
 =&\{K,  T_{0}, T_{0}\}+\{K,  T_{0}, \oplus_{\alpha \in \Lambda^{\neg J}} T_{\alpha}\}+\{K,  T_{0}, \oplus_{\beta \in \Lambda^{J}} T_{\beta} \}\\
 +&\{K,  \oplus_{\alpha \in \Lambda^{\neg J}} T_{\alpha}, T_{0}\}+\{K,  \oplus_{\alpha \in \Lambda^{\neg J}} T_{\alpha}, \oplus_{\gamma \in \Lambda^{\neg J}} T_{\gamma}\}+\{K,  \oplus_{\alpha \in \Lambda^{\neg J}} T_{\alpha}, \oplus_{\beta \in \Lambda^{J}} T_{\beta}\}\\
 +&\{K,  \oplus_{\beta \in \Lambda^{ J}} T_{\beta}, T_{0}\}+\{K,  \oplus_{\beta \in \Lambda^{ J}} T_{\beta}, \oplus_{\alpha \in \Lambda^{\neg J}} T_{\alpha}\}+\{K,  \oplus_{\beta \in \Lambda^{ J}} T_{\beta}, \oplus_{\gamma \in \Lambda^{J}} T_{\gamma}\}.
 \end{align*}
Here, it is easy to see $$\{K,  T_{0}, \oplus_{\beta \in \Lambda^{J}} T_{\beta} \}=0,$$
$$\{K,  \oplus_{\alpha \in \Lambda^{\neg J}} T_{\alpha}, \oplus_{\beta \in \Lambda^{J}} T_{\beta}\}=0,$$
$$\{K,  \oplus_{\beta \in \Lambda^{ J}} T_{\beta}, T_{0}\}=0,$$
$$\{K,  \oplus_{\beta \in \Lambda^{ J}} T_{\beta}, \oplus_{\alpha \in \Lambda^{\neg J}} T_{\alpha}\}=0,$$
$$\{K,  \oplus_{\beta \in \Lambda^{ J}} T_{\beta}, \oplus_{\gamma \in \Lambda^{J}} T_{\gamma}\}=0.$$
So
\begin{align*}
 &\{K,T,T\}\\
 =&\{K,  T_{0}\oplus(\oplus_{\alpha \in \Lambda^{\neg J}} T_{\alpha})\oplus(\oplus_{\beta \in \Lambda^{J}} T_{\beta}), T_{0}\oplus(\oplus_{\alpha \in \Lambda^{\neg J}} T_{\alpha})\oplus(\oplus_{\beta \in \Lambda^{J}} T_{\beta})\}\\
 =&\{K,  T_{0}, T_{0}\}+\{K,  T_{0}, \oplus_{\alpha \in \Lambda^{\neg J}} T_{\alpha}\}+\{K,  \oplus_{\alpha \in \Lambda^{\neg J}} T_{\alpha},T_{0}\}
 +\{K,  \oplus_{\alpha \in \Lambda^{\neg J}} T_{\alpha}, \oplus_{\gamma \in \Lambda^{\neg J}} T_{\gamma}\}.
  \end{align*}

It is easy to get $\{K,  T_{0}, T_{0}\}\subset K$. Let us consider $\{K,  T_{0}, \oplus_{\alpha \in \Lambda^{\neg J}} T_{\alpha}\}$ and suppose there exist $\alpha_{i}
\in  \Lambda^{ J,I}$ and $\alpha \in \Lambda^{ \neg J}$ such that $\{T_{-\alpha_{i}}, T_{0}, T_{\alpha}\}\neq 0$. Since $T_{-\alpha_{i}}\subset K \subset J$, we get
$-\alpha_{i}+\alpha \in \Lambda^{ J}$. By the root-multiplicativity of $T$,
the symmetries of $\Lambda^{ J}$ and $\Lambda^{\neg J}$, and the fact $T_{\alpha_{i}}\subset I$ we obtain $0 \neq \{T_{\alpha_{i}}, T_{0}, T_{-\alpha}\}=T_{\alpha_{i}-\alpha}\subset I$, that is $\alpha_{i}-\alpha \in \Lambda^{ J,I}$. Hence, $-\alpha_{i}+\alpha \in -\Lambda^{ J,I}$ and so
$\{T_{-\alpha_{i}}, T_{0}, T_{\alpha}\}=T_{-\alpha_{i}+\alpha}\subset K$. Similarly, we also get $\{K,   \oplus_{\alpha \in \Lambda^{\neg J}}, T_{\alpha}, T_{0}\}\subset K.$

At last, we consider $\{K,  \oplus_{\alpha \in \Lambda^{\neg J}} T_{\alpha}, \oplus_{\gamma \in \Lambda^{\neg J}} T_{\gamma}\}$ and suppose there exist
$\alpha_{i}
\in  \Lambda^{ J,I},$  $\alpha \in \Lambda^{ \neg J}$ and $\gamma \in \Lambda^{ \neg J}$ such that $\{T_{-\alpha_{i}}, T_{\alpha}, T_{\gamma}\}\neq 0.$ Since $T_{-\alpha_{i}}\subset K \subset J$, we get
$-\alpha_{i}+\alpha+\gamma \in \Lambda^{ J}$.   By the root-multiplicativity of $T$,
the symmetries of $\Lambda^{ J}$ and $\Lambda^{\neg J}$, and the fact $T_{\alpha_{i}}\subset I$, one gets  $0 \neq \{T_{\alpha_{i}}, T_{-\alpha}, T_{-\gamma}\}=T_{\alpha_{i}-\alpha-\gamma}\subset I$, that is $\alpha_{i}-\alpha-\gamma \in \Lambda^{ J,I}$.  Hence, $-\alpha_{i}+\alpha+\gamma \in -\Lambda^{ J,I}$
and so
$\{T_{-\alpha_{i}}, T_{\alpha}, T_{\gamma}\}=T_{-\alpha_{i}+\alpha+\gamma}\subset K$. Consequently, $K$ is an ideal of $T$.
\epf

We introduce the Definition of primeness in the framework of Leibniz triple systems following the same
motivation that in the case of simplicity (see Definition \ref{ 3.16789}  and the above paragraph).

\bdefn
A Leibniz triple system $T$ is said to be prime if given two ideals $I$, $K$ of $T$ satisfying $\{I, K,I\} +
\{K, I,I\}+\{I, I, K\} = 0$, then either $I \in \{0, J, T\}$ or $K \in \{0, J, T\}$.
\edefn

We also note that the above Definition agrees with the Definition of prime Lie triple system, since $J= {0}$
in this case.

Under the hypotheses of Proposition \ref{67891110345777777999} we have:

\bcor\label{555}
If furthermore $T$ is prime, then any nonzero ideal $I$ of $T$ such that $I\subseteq J$ satisfies $I= J$.
\ecor
\bpf
Observe that, by Proposition \ref{67891110345777777999}, we could have $J = I \oplus K$ with $I, K$  ideals of $T$, being $\{I, K, I\} +
\{K, I, I\}+ \{ I, I, K\} = 0$ as consequence of $I, K \subseteq J$. The primeness of $T$ completes the proof.
\epf

\bprop \label{67877999}
Suppose $T= \{T, T,T\}$, $\mathrm{Ann}_{\mathrm{Lie}} (T)=0$, $T$ is root-multiplicative and for any $\alpha\in \Lambda^{1}$, we have $ \mathrm{dim}L_{\alpha}^{0}=1 $. If $\Lambda^{\neg J}$
 has all of its roots $\neg J$-connected and $\{T_{0}, T_{\alpha}, T_{\beta}\}=\{T_{\alpha}, T_{0},  T_{\beta}\}=\{T_{\alpha}, T_{\beta}, T_{0}\}=0,$ for $\alpha, \beta \in  \Lambda^{\neg J}$,
then any ideal $I$ of $T$ such that $I \not \subseteq  J$ satisfies $I = T$.
\eprop
\bpf
Taking into account Lemma \ref{lemma 4.4} and Proposition \ref{67890345777777999} we just have to study the case in which
$$I=(I\cap T_{0})\oplus(\oplus_{\beta_{j} \in \Lambda^{ J,I}} T_{\beta_{j}}),$$
\noindent where $I\cap T_{0}\neq 0$. But this possibility never happens. Indeed, observe that
$$\{I\cap T_{0}, T_{0}, T_{0}\}+\{T_{0}, I\cap T_{0},  T_{0}\}+\{ T_{0}, T_{0}, I\cap T_{0}\}=0.$$
We also have $$\{I\cap T_{0}, \oplus_{\alpha \in \Lambda^{\neg J}} T_{\alpha}, T_{0}\}+\{ \oplus_{\alpha \in \Lambda^{\neg J}} T_{\alpha}, I\cap T_{0}, T_{0}\}
+\{ \oplus_{\alpha \in \Lambda^{\neg J}} T_{\alpha},  T_{0}, I\cap T_{0}\}\subset I\cap \oplus_{\alpha \in \Lambda^{\neg J}} T_{\alpha}=0$$ and  $$\{I\cap T_{0}, T_{0}, \oplus_{\alpha \in \Lambda^{\neg J}} T_{\alpha}\}+\{ T_{0}, I\cap T_{0}, \oplus_{\alpha \in \Lambda^{\neg J}} T_{\alpha} \}
+\{ T_{0}, \oplus_{\alpha \in \Lambda^{\neg J}} T_{\alpha},   I\cap T_{0}\}\subset I\cap \oplus_{\alpha \in \Lambda^{\neg J}} T_{\alpha}=0.$$
And known condition gives us $$\{I\cap T_{0}, \oplus_{\alpha \in \Lambda^{\neg J}} T_{\alpha}, \oplus_{\beta \in \Lambda^{\neg J}} T_{\beta}\}+\{ \oplus_{\alpha \in \Lambda^{\neg J}} T_{\alpha}, I\cap T_{0}, \oplus_{\beta \in \Lambda^{\neg J}} T_{\beta}\}
+\{ \oplus_{\alpha \in \Lambda^{\neg J}} T_{\alpha},  \oplus_{\beta \in \Lambda^{\neg J}} T_{\beta}, I\cap T_{0}\}$$
=0.

That is, we get $I\cap T_{0}\subset \mathrm{Ann}_{\mathrm{Lie}} (T)=0$, a contradiction. Proposition \ref{67890345777777999} completes the proof.
\epf

Given any $\alpha \in \Lambda^{\gamma}$, $\gamma \in \{J, \neg J\}$, we denote by
$$\Lambda_{\alpha}^{\gamma}:= \{\beta \in  \Lambda^{\gamma}: \beta \sim_{\neg J} \alpha\}.$$
For $\alpha \in  \Lambda^{\gamma}$, we write $T_{0,\Lambda_{\alpha}^{\gamma}}:=\mathrm{span}_{\mathbb{K}}\{\{T_{\tau}, T_{\beta}, T_{\gamma}\}: \tau+\beta+\gamma=0;
 \tau,\beta,\gamma \in \Lambda_{\alpha}^{\gamma}\}\subset T_{0},$ and $V_{\Lambda_{\alpha}^{\gamma}}:=\oplus_{\beta \in \Lambda_{\alpha}^{\gamma}}T_{\beta}$. We denote by $T_{\Lambda_{\alpha}^{\gamma}}$ the following subspace of $T$, $T_{\Lambda_{\alpha}^{\gamma}}:=T_{0,\Lambda_{\alpha}^{\gamma}}\oplus V_{\Lambda_{\alpha}^{\gamma}}.$

\blem\label{2222}
If $T=\{T, T, T\}$, then $T_{\Lambda_{\alpha}^{J}}$ is an ideal of $T$ for any $\alpha \in \Lambda^{J}$.
\elem

\bpf
We have  $T_{0,\Lambda_{\alpha}^{J}}=0$ and so
$$T_{\Lambda_{\alpha}^{J}}=\oplus_{\beta \in \Lambda_{\alpha}^{J}}T_{\beta}.$$
Using Proposition \ref{38888} (1)  and (\ref{8889765566897756}), we have
\begin{align*}
 &\{T_{\Lambda_{\alpha}^{J}},T,T\}+\{T,T_{\Lambda_{\alpha}^{J}},T\}+\{T,T,T_{\Lambda_{\alpha}^{J}}\}\\
=&\{T_{\Lambda_{\alpha}^{J}},T,T\}+0+0\\
=&\{T_{\Lambda_{\alpha}^{J}},T,T\}\\
=&\{T_{\Lambda_{\alpha}^{J}},  T_{0}\oplus(\oplus_{\alpha \in \Lambda^{\neg J}} T_{\alpha})\oplus(\oplus_{\beta \in \Lambda^{J}} T_{\beta}), T_{0}\oplus(\oplus_{\alpha \in \Lambda^{\neg J}} T_{\alpha})\oplus(\oplus_{\beta \in \Lambda^{J}} T_{\beta})\}\\
 =&\{T_{\Lambda_{\alpha}^{J}},  T_{0}, T_{0}\}+\{ T_{\Lambda_{\alpha}^{J}}, \oplus_{\gamma \in \Lambda^{\neg J}}T_{\gamma}, T_{0}\}+\{ T_{\Lambda_{\alpha}^{J}}, T_{0}, \oplus_{\gamma \in \Lambda^{\neg J}}T_{\gamma}\}+\{T_{\Lambda_{\alpha}^{J}},\oplus_{\beta\in \Lambda^{\neg J}}T_{\beta}, \oplus_{\gamma \in \Lambda^{\neg J}}T_{\gamma}\}.
 \end{align*}
It is easy to see $\{T_{\Lambda_{\alpha}^{J}}, T_{0}, T_{0}\}\subset T_{\Lambda_{\alpha}^{J}}.$ Next we will assert $\{ T_{\Lambda_{\alpha}^{J}}, \oplus_{\gamma \in \Lambda^{\neg J}}T_{\gamma}, T_{0}\}\subset T_{\Lambda_{\alpha}^{J}}$. Indeed, given any $m \in \Lambda_{\alpha}^{J}$, $\gamma \in \Lambda^{\neg J}$ such that $\{T_{m}, T_{\gamma},  T_{0}\}\neq 0$. By $m+\gamma\neq 0$, we have $m+\gamma \in \Lambda^{J}$ and  so $\{m, \gamma, 0\}$ is a $\neg J$-connection from $m$ to $m+\gamma$. By the symmetry and transitivity of $\sim_{\neg J}$ in $\Lambda^{J}$,   we get $m+\gamma \in \Lambda_{\alpha}^{J}.$ Hence $\{ T_{\Lambda_{\alpha}^{J}}, \oplus_{\gamma \in \Lambda^{\neg J}}T_{\gamma}, T_{0}\}\subset T_{\Lambda_{\alpha}^{J}}$. Similarly, we also get $\{ T_{\Lambda_{\alpha}^{J}}, T_{0}, \oplus_{\gamma \in \Lambda^{\neg J}}T_{\gamma}\}\subset T_{\Lambda_{\alpha}^{J}}$.

 At last, we will assert $\{T_{\Lambda_{\alpha}^{J}},\oplus_{\beta\in \Lambda^{\neg J}}T_{\beta}, \oplus_{\gamma \in \Lambda^{\neg J}}T_{\gamma}\}\subset T_{\Lambda_{\alpha}^{J}}$. Indeed, given any $m \in \Lambda_{\alpha}^{J}$, $\beta, \gamma \in \Lambda^{\neg J}$ such that
$\{T_{m}, T_{\beta},  T_{\gamma}\}\neq 0$. By $m+\beta+\gamma\neq 0$, we have $m+\beta+\gamma \in \Lambda^{J}$ and so $\{m, \beta, \gamma\}$ is a $\neg J$-connection from $m$ to $m+\beta+\gamma$.  By the symmetry and transitivity of $\sim_{\neg J}$ in $\Lambda^{J}$,    we get $m+\beta+\gamma \in \Lambda_{\alpha}^{J}.$ Hence $\{T_{\Lambda_{\alpha}^{J}},\oplus_{\beta\in \Lambda^{\neg J}}T_{\beta}, \oplus_{\gamma \in \Lambda^{\neg J}}T_{\gamma}\}\subset T_{\Lambda_{\alpha}^{J}}$.

Consequently, $\{T_{\Lambda_{\alpha}^{J}},T,T\}+\{T,T_{\Lambda_{\alpha}^{J}},T\}+\{T,T,T_{\Lambda_{\alpha}^{J}}\}\subset  T_{\Lambda_{\alpha}^{J}}$. So $T_{\Lambda_{\alpha}^{J}}$ is an ideal of $T$ for any $\alpha \in \Lambda^{J}$.
\epf

\bthm
Suppose $T = \{T, T, T\}$, $\mathrm{Ann}_{\mathrm{Lie}}(T) = 0$, $T$ is root-multiplicative and for any $\alpha\in \Lambda^{1}$,  we have $ \mathrm{dim}L_{\alpha}^{0}=1 $. If $\Lambda^{J}$, $\Lambda^{\neg J}$ are symmetric and $\{T_{0}, T_{\alpha}, T_{\beta}\}=\{T_{\alpha}, T_{0},  T_{\beta}\}=\{T_{\alpha}, T_{\beta}, T_{0}\}=0,$ for $\alpha, \beta \in  \Lambda^{\neg J}$,
then $T$ is simple if and only if it is prime and $\Lambda^{J}$, $\Lambda^{\neg J}$ have all of their roots $\neg J$-connected.
\ethm

\bpf
Suppose $T$ is simple. If $\Lambda^{J}\neq \emptyset$ and we take $\alpha \in \Lambda^{J}$.  Lemma \ref{2222} gives us $T_{\Lambda_{\alpha}^{J}}$
is a nonzero ideal
of $T$ and so, (by simplicity), $T_{\Lambda_{\alpha}^{J}}
= J = \oplus_{\beta \in \Lambda^{J}} T_{\beta}$ (see  (\ref{888976556689}) and (\ref{888778555})). Hence, $\Lambda_{\alpha}^{J}=\Lambda^{J}$
 and consequently $\Lambda^{J}$ has all of its roots $\neg J$-connected.

Consider now any $\gamma \in \Lambda^{\neg J}$ and the subspace $T_{\Lambda_{\gamma}^{\neg J}}$.
Let us denote by $I(T_{\Lambda_{\gamma}^{\neg J}})$ the ideal of $T$ generated by $T_{\Lambda_{\gamma}^{\neg J}}$.
We observe that the fact $J$ is an ideal of $T$ and we assert that $I(T_{\Lambda_{\gamma}^{\neg J}})\cap (\oplus_{\delta \in \Lambda^{\neg J}}T_{\delta})$ is contained
in the linear span of the set
\begin{align*}
 &\Big\{\{\{\cdots\{v_{\gamma^{'}},v_{\alpha_{1}}, v_{\alpha_{2}}\},\cdots\}, v_{\alpha_{2n}},v_{\alpha_{2n+1}}\}; \{v_{\alpha_{2n+1}},  v_{\alpha_{2n}},\{\cdots\{ v_{\alpha_{2}},v_{\alpha_{1}}, v_{\gamma^{'}}\},\cdots\}\};\\
 &\{\{\cdots\{v_{\alpha_{1}}, v_{\alpha_{2}},v_{\gamma^{'}}\},\cdots\}, v_{\alpha_{2n}},v_{\alpha_{2n+1}}\}; \{v_{\alpha_{2n+1}},  v_{\alpha_{2n}},\{\cdots\{ v_{\gamma^{'}}, v_{\alpha_{2}},v_{\alpha_{1}}\},\cdots\}\};\\
 &\{\{\cdots\{,v_{\alpha_{1}},v_{\gamma^{'}}, v_{\alpha_{2}}\},\cdots\}, v_{\alpha_{2n}},v_{\alpha_{2n+1}}\}; \{v_{\alpha_{2n+1}},  v_{\alpha_{2n}},\{\cdots\{ v_{\alpha_{2}},v_{\gamma^{'}},v_{\alpha_{1}}, \},\cdots\}\};\\
 & \mathrm{where} \ 0\neq v_{\gamma^{'}} \in T_{\Lambda_{\gamma}^{\neg J}}, 0\neq v_{\alpha_{i}}\in T_{\alpha_{i}}, \alpha_{i} \in \Lambda^{\neg J} \mathrm{and}\  n \in \mathbb{N}\Big\}.
 \end{align*}

By simplicity $I(T_{\Lambda_{\gamma}^{\neg J}})=T$. From here, given any $\delta \in \Lambda^{\neg J}$, the above observation gives us that we can write $\delta=\gamma^{'}+\alpha_{1}+\cdots+\alpha_{2n+1}$ for $\gamma^{'} \in \Lambda_{\gamma}^{\neg J}$, $\alpha_{i} \in \Lambda^{\neg J}$ and being the partial sums nonzero. Hence $\{\gamma^{'},\alpha_{1},\cdots,\alpha_{2n+1}\}$ is a $\neg J$-connection from $\gamma^{'}$ to $\delta$. By the symmetry and transitivity of $\sim_{\neg J}$ in $\Lambda^{\neg J}$, we deduce $\gamma$ is $\neg J$-connected to any $\delta \in  \Lambda^{\neg J}$. Consequently, Proposition \ref{67890345777777} gives us that
$\Lambda^{\neg J}$ has all of its roots $\neg J$-connected.
Finally, since $T$ is simple then is prime.

The converse is consequence of Corollary \ref{555}  and Proposition \ref
{67877999}.
\epf

\end{document}